\documentclass[12pt]{amsart}
\usepackage{amssymb}
\newtheorem*{theo}{Theorem}
\newtheorem*{cor}{Corollary}
\newtheorem*{prop}{Proposition}
\newtheorem*{lem}{Lemma}
\newtheorem*{poslem}{Positivity lemma}
\theoremstyle{remark}
\newtheorem*{rem}{Remark}

\DeclareMathOperator{\gr}{gr}
\DeclareMathOperator{\car}{Char}
\DeclareMathOperator{\SP}{SP}
\DeclareMathOperator{\Perv}{Perv}
\DeclareMathOperator{\graph}{graph}
\renewcommand{\ker}{\mathop{\rm Ker}\nolimits}
\DeclareMathOperator{\id}{Id}

\def\pQQ{{{}^p\!\QQ}}
\def\pCC{{{}^p\!\CC}}
\def\psip{{{}^p\!\psi}}
\def\phip{{{}^p\!\phi}}

\def\cFG{\ccF_{\!G}}
\def\FGir{\cFG^{\rm ir}}

\def\per{{}^p\!}

\def\crit{\mathop{\rm Crit}\nolimits}

\newcommand{\ootimes}{\mathop\otimes\limits}
\newcommand{\lefpar}{\left(}
\newcommand{\rigpar}{\right)}

\newcommand{\module}[1]{\left\vert#1\right\vert}
\newcommand{\norme}[1]{\left\Vert#1\right\Vert}
\let\dpl\displaystyle

\let\wt\widetilde
\let\ov\overline

\newcommand{\defin}{:=}

\newcommand{\loccit}{{\em loc.\ cit}}
\newcommand{\eg}{{\it e.g}}

\newcommand{\ie}{{\it i.e}}
\newcommand{\resp}{{\it resp}}

\newcommand{\T}{\S\kern .15em }
\newcommand{\ptbl}{.\kern .15em }
\newcommand{\bbullet}{{\scriptscriptstyle\bullet}}

\DeclareMathAlphabet{\mathcalmaigre}{U}{eus}{m}{n}

\def\AA{\mathbf{A}}
\def\CC{\mathbb{C}}
\def\ZZ{\mathbb{Z}}
\def\QQ{\mathbb{Q}}
\def\NN{\mathbb{N}}
\def\RR{\mathbb{R}}
\def\PP{\mathbb{P}}

\def\bR{\boldsymbol{R}}

\def\cD{\mathcal{D}}

\def\cH{\mathcal{H}}

\def\cM{\mathcal{M}}
\def\cN{\mathcal{N}}
\def\cO{\mathcal{O}}

\def\cU{\mathcal{U}}

\def\cW{\mathcal{W}}
\def\cX{\mathcalmaigre{X}}
\def\cY{\mathcal{Y}}

\def\ccF{\mathcalmaigre{F}}
\def\ccG{\mathcalmaigre{G}}

\def\Afuc{\check\AA\!^1}

\newdimen\lengtharrow
 \newbox\exponantbox \newbox\indicebox
\def\dimmax#1#2{\ifdim#1<#2 #2\else #1\fi}

\def\arrowr#1#2%
    {\setbox\exponantbox=\hbox{$#2$}
      \setbox\indicebox=\hbox{$#1$}
       \lengtharrow=\dimmax{\wd\indicebox}{\wd\exponantbox}
        \ifdim \lengtharrow=0pt
															\mathrel{\smash{\mathop{\hbox to 10mm{\rightarrowfill}}}}
         \else \ifdim\lengtharrow<6mm
                \lengtharrow=10mm
                 \else \advance\lengtharrow by 4mm
                  \fi
					\mathrel{\smash{\mathop{\hbox to\lengtharrow{\rightarrowfill}}\limits%
                         ^{\textstyle#2}_{\textstyle#1}}}
									\fi}

\def\arrowl#1#2%
    {\setbox\exponantbox=\hbox{$#2$}
      \setbox\indicebox=\hbox{$#1$}
       \lengtharrow=\dimmax{\wd\indicebox}{\wd\exponantbox}
        \ifdim \lengtharrow=0pt
															\mathrel{\smash{\mathop{\hbox to 10mm{\leftarrowfill}}}}
         \else \ifdim\lengtharrow<6mm
                \lengtharrow=10mm
                 \else \advance\lengtharrow by 4mm
                  \fi
					\mathrel{\smash{\mathop{\hbox to\lengtharrow{\leftarrowfill}}\limits%
                         ^{\textstyle#2}_{\textstyle#1}}}
									\fi}

\def\arrowd#1#2{\llap{$ \textstyle #1$}\left\downarrow%
\vbox to 6mm{}\right.\rlap{$\textstyle #2$}}

\def\MRE#1{\arrowr{}{#1}}

\def\MDR#1{\arrowd{}{#1}}
\def\MDL#1{\arrowd{#1}{}}

\newcommand{\isom}{\stackrel{\sim}{\longrightarrow}}

\begin{document}

\title{Semicontinuity of the spectrum at infinity}
\author{Andr\'as N\'emethi}
\address{Department of Mathematics, the Ohio State University, 231 West, 18th Avenue, Columbus OH 43210, U.S.A.}
\email{nemethi@math.ohio-state.edu}
\author{Claude Sabbah}
\address{URA 169 du C.N.R.S., Centre de Math\'ematiques, \'Ecole polytechnique, F--91128 Palaiseau cedex, France}
\email{sabbah@math.polytechnique.fr}
\begin{abstract}
We prove that, for an analytic family of ``weakly tame'' regular functions on an affine manifold, the spectrum at infinity of each function of the family is semicontinuous in the sense of Varchenko.
\end{abstract}
\date{July 1997}
\maketitle

\section*{Introduction}
Let $U$ be an affine manifold and let $f:U\rightarrow \CC$ be a nonconstant regular function with only isolated critical points in $U$.

We say that $f$ is {\em cohomologically tame} if there exists an extension $\ov f:X\rightarrow \CC$ of $f$ with $X$ quasiprojective and $\ov f$ proper, such that for any $c\in\CC$, the support of the vanishing cycle sheaf of the function $\ov f-c$ with coefficients in the sheaf $j_!\QQ_U$, or equivalently in the complex $\bR j_*\QQ_U$, does not meet $X-U$ ($j$ is the inclusion of $U$ in $X$ and $j_!$ denotes the extension by $0$). In particular $f$ is onto. This condition is satisfied if and only if $j_!\QQ_U$ or equivalently $\bR j_*\pQQ_U$ is noncharacteristic with respect to $\ov f$ (see \eg.\ \cite[prop-def\ptbl1.1]{Parusinski95}), \ie., choosing an embedding of $\ov f$ into $F:\cX\rightarrow \CC$ with $\cX$ smooth, there are no points $(x,dF(x))$ in the characteristic variety (or microssuport, see \cite{K-S90}) $\car j_!\QQ_U\subset T^*\cX$ such that $x\in X-U$ and $dF(x)\neq0$.

We say that $f$ is {\em M-tame} (\cite{N-Z90}) if for some closed embedding $U\subset\CC^N$ (\ie.\ for some presentation of $\cO(U)$ as a finite type $\CC$-algebra) and some $a\in\CC^N$, if $\delta$ denotes the distance function $\delta(x)=\norme{x-a}^2$, then for any $\eta>0$ there exists $R(\eta)>0$ such that, for any $r\geq R(\eta)$, the spheres $\delta(x)=r$ are transversal to $f_{}^{-1}(t)$ for $\module{t}\leq \eta$.

M-tameness is a property analogous, in the affine case, to the existence of the Milnor fibration of hypersurface singularities. It emphasizes the behaviour of $f$ on the affine manifold $U$. Cohomological tameness is more algebraic and emphasizes the behaviour of $f$ at infinity for some compactification: there is only one ``Milnor ball'', namely $U$ itself. The latter notion is only cohomological (with $\QQ$ as coefficients) but the former is topological.

It is not clear whether one property is stronger than the other one. It is known that for polynomials on $\CC_{}^{n+1}$, cohomological tameness with respect to the standard projective compactification of the fibers is equivalent to the so called Malgrange condition (\cite{Parusinski95}), which in turn implies M-tameness. For instance, tame polynomials on $\CC_{}^{n+1}$ in the sense of Broughton \cite{Broughton88} are both M-tame and cohomologically tame. Conversely, there exist polynomials which are both M-tame and cohomologically tame (with respect to some compactification) but do not satisfy the Malgrange condition (\ie.\ are not cohomologically tame with respect to the standard projective compactification): many examples of such polynomials have been constructed by L. P\u{a}unescu and A. Zaharia.

In the following we will call {\em weakly tame} a function which is either cohomologically tame or M-tame. Let us mention, as another example of a weakly tame function, the function induced on the affine smooth manifold $U\subset\CC^N$ by a linear pencil on $\PP^N$ which basis $L\subset\PP_{\infty}^{N-1}$ is transversal to a Whitney stratification of the closure of $U$ in $\PP^N$ containing $U$ as a stratum.

\medskip
Denote $\mu(f)$ the sum of the Milnor numbers of $f$ at its critical points. It can be decomposed as $\mu(f)=\sum_{\beta\in\QQ}^{}\nu_\beta(f)$, where $\beta\mapsto\nu_\beta(f)$ denotes the integral valued function associated with {\em the spectrum of $f$ at infinity} (see \cite{Bibi96b} for cohomologically tame functions and see below \T\ref{secMtame} for M-tame functions). It is known that this function satisfies $\nu_\beta=0$ for $\beta\not\in[0,\dim U]$ (and for $\beta\not\in{}]0,\dim U[$ if $U=\CC_{}^{n+1}$) and $\nu_{\dim U-\beta}^{}=\nu_{\beta}^{}$ (see \loccit.). For any $\gamma\in\QQ$ we denote $\Sigma_\gamma(f)=\dpl\sum_{\beta\in{}]\gamma,\gamma+1]}^{}\nu_\beta(f)$.

\medskip
Let now $(S,0)$ be the analytic germ of a smooth curve, let $\pi:\cU\rightarrow S$ be a smooth affine morphism with smooth affine fibres $\cU_s\defin\pi_{}^{-1}(s)$, and let $f:\cU\rightarrow \CC$ be an analytic family of regular functions on $\cU_s$. Assume that $f_s:\cU_s\rightarrow \CC$ is weakly tame for any $s\in S$. Then one has $\mu(f_0)\leq \mu(f_s)$ ($s\in S$): indeed, let $C$ be the critical locus of $(f,\pi):\cU\rightarrow \CC\times S$; the assumption implies that $C$ is locally a complete intersection curve and $(f,\pi)$ is quasi-finite on it; thus $\mu(f_0)=\mu'(f_s)$ where $\mu'$ denotes the sum over the critical points lying on the components $C'$ of $C$ which intersect $\cU_0$ (and on which $(f,\pi)$ is finite).

The purpose of this note is to prove an analogue of the semicontinuity of the spectrum of an isolated hypersurface singularity, a theorem proved by Varchenko \cite{Varchenko83} and Steenbrink \cite{Steenbrink85}, namely

\begin{theo}
Let $f:\cU\rightarrow \CC$ be an analytic family of regular functions on $\cU_s$. Assume that for any $s\in S$ the function $f_s:\cU_s\rightarrow \CC$ is weakly tame. Then for any $\gamma\in\QQ$ we have $\Sigma_\gamma(f_0)\leq \Sigma_\gamma(f_s)$.
\end{theo}

Throughout this note we use notation, conventions and results of \cite{Bibi96b,Bibi96a}. In particular, we refer to the appendix of \cite{Bibi96a} for the convention concerning perverse functors.

The proof is given in \T\ref{secproofth} and follows the ideas in \cite{Varchenko83} and \cite{Steenbrink85}. In \cite{G-N96} a partial result for $(*)$-polynomials (defined in \cite{G-N95}) is obtained by a method which reduces to the theorem of Varchenko and Steenbrink.

\medskip
The first author thanks A. Zaharia for helpful discussions and the second author thanks J.H.M. Steenbrink for raising this problem to his attention.

\section{M-tame functions}\label{secMtame}
We indicate in this section why the results of \cite{Bibi96b} for cohomologically tame functions apply as well to M-tame functions. The main points are theorems 8.1, 10.1 and corollary 11.1 in \cite{Bibi96b}, from which follow all the results in \loccit. \T13. In the following we assume that $f$ is M-tame.

\medskip
The main ingredient in the proof of theorem 8.1 in \loccit. is the fact that $\bR f_*$ commutes with $\psi_{f-c}^{}$ and $\phi_{f-c}^{}$ on $\QQ_U$ for any $c\in\CC$ (see \loccit. cor\ptbl 8.4). In the case of M-tame functions this follows from the fact that for any disc $D$ in $\CC$, there exists $R_0>0$ such that for $R>R_0$, putting $B(R)=U\cap\{\delta(x)\leq r\}$, the restriction morphism $\bR f_*\QQ_{f_{}^{-1}(D)}^{}\rightarrow \bR f_*\QQ_{B_R\cap f_{}^{-1}(D)}^{}$ is an isomorphism (which follows easily from the definition). Indeed, as the map $f:B_R\cap f_{}^{-1}(D)\rightarrow D$ is proper, one may apply to it the commutation of $\psi_{f-c}^{},\phi_{f-c}^{}$ with $\bR f_*$.

\medskip
Before proving the coherency of the Brieskorn lattice, let us prove that $f_+(\cO_{U}^{\rm an})$ has regular holonomic cohomology and that the natural morphism $(f_+\cO_U)_{}^{\rm an}\rightarrow f_+(\cO_{U}^{\rm an})$ is an isomorphism.

The coherency of the cohomology of $f_+(\cO_{U}^{\rm an})$ as a $\cD_{\CC}^{\rm an}$-module is a direct corollary of \cite[cor\ptbl 8.1]{Sch-Sch94} by taking there for $S$ a point, for $Y$ a disc, $f_{}^{-1}(Y)$ for $X$, $\cO_{X}^{\rm an}$ for $\cM$, the function $\delta$ for $\varphi$ and the constant sheaf $\CC_U$ for $F$. Assumptions (i) and (ii) of \loccit. are clearly satisfied and (iii) is the M-tameness condition. Kashiwara's estimate for the characteristic variety shows that the cohomology is indeed holonomic.

\medskip
Let us quickly recall the regularity of $f_+(\cO_{U}^{\rm an})$. One has to show that for any $c\in \CC$, if $t$ is a local coordinate centered at $c$, there exits a $\CC\{t\}\langle t\partial_t\rangle$-submodule of the germ of $f_+(\cO_{U}^{\rm an})$ at $c$, generating it as a $\CC\{t\}\langle\partial_t\rangle$-module, and which is finitely generated as a $\CC\{t\}$-module.

Let $i_f:U\hookrightarrow U\times \CC$ be the inclusion of the graph and let $p$ be the second projection. Consider a coherent $\cD_{U\times \CC/\CC}^{\rm an}$-module $\cN_{}^{\rm an}$ generating $i_{f+}^{}\cO_{U}^{\rm an}$ as a $\cD_{U\times \CC}^{\rm an}$-module. We may now apply the same result with $S=Y$ a disc, $X=U\times \CC$, $F=\CC_X$ and $\cM=\cN_{}^{\rm an}$, $\varphi=\delta$: outside of the critical points of $f$, the characteristic variety of $\cN_{}^{\rm an}$ is the relative conormal space to the map $p$ restricted to the graph of $f$ and the tameness condition gives (iii) in \loccit. It follows that $p_+\cN_{}^{\rm an}$ has $\cO_{\CC}^{\rm an}$-coherent cohomology. 

As $i_{f+}^{}\cO_{U}^{\rm an}$ is regular holonomic, there exists such an $\cN_{}^{\rm an}$ which is stable by the action of $t\partial_t$ (see \eg.\ \cite{M-S86}). The image of $\cH^0(p_+\cN_{}^{\rm an})$ in $\cH^0(f_+(\cO_{U}^{\rm an}))$ will be the desired $\CC\{t\}\langle t\partial_t\rangle$-submodule.

\medskip
Let now $\ov f: X\rightarrow \CC$ be a smooth quasi-projective compactification of $f$ and let $j:U\hookrightarrow X$ be the inclusion. One has $f_+\cO_U=\ov f_+(j_+\cO_U)$, $(f_+\cO_U)_{}^{\rm an}=\ov f_+(j_+\cO_U)_{}^{\rm an}$ and $f_+(\cO_{U}^{\rm an})=\ov f_+(\bR j_*\cO_U)_{}^{\rm an}$. The natural morphism $\ov f_+(j_+\cO_U)_{}^{\rm an}\rightarrow \ov f_+(\bR j_*\cO_U)_{}^{\rm an}$ induces an isomorphism of the corresponding de~Rham complexes as the analytic de~Rham functor commutes with $\ov f_+$ and $\bR f_*$ (see \eg.\ \cite[II.5.5]{M-N90}) and because the holonomic $\cD_X$-module $j_+\cO_U$ is regular by the Grothendieck comparison theorem. This natural morphism is then an isomorphism according to the Riemann-Hilbert correspondence in dimension $1$.

\medskip
Let us now prove the coherency of the Brieskorn lattice $M_0$. Remark first that as $U_{}^{\rm an}$ is Stein, the object $f_+\cO_{U}^{\rm an}$ is  the complex $(\Omega_{}^{\bbullet+\dim U}(U_{}^{\rm an})[\partial _t],d_f)$, where $\Omega^k(U_{}^{\rm an})$ is the space of holomorphic $k$-forms on $U_{}^{\rm an}$ and $d_f$ is the twisted differential $d_f\omega=d\omega-df\wedge\omega\partial _t$. In order to prove the analogue of \cite[th\ptbl 10.1]{Bibi96b}, one first shows the $\cO_{\CC}^{\rm an}$-coherency of the image $\cM_0$ of $\Omega_{}^{\dim U}(U_{}^{\rm an})$ in $\cH^0f_+\cO_{U}^{\rm an}$. This can be done as in \loccit. using the coherency for the direct image of relative submodules used above.

Let then $M_0$ be the image of $\Omega_{}^{\dim U}(U)$ in $\cH^0f_+\cO_U$. We have $M_0\subset\Gamma(\CC,\cM_0)\cap \cH^0f_+\cO_U$ and the latter term has finite type over $\CC[t]$, hence so does $M_0$.

\medskip
The proof of the duality theorem 11.1 of \loccit. remains valid in this context. \hfill\qed

\section{Thom-Sebastiani theorem for the spectrum at infinity} \label{secTS}

\begin{prop}
If $f:U\rightarrow \CC$, $g:V\rightarrow \CC$ are both M-tame (\resp.\ both cohomologically tame) functions on affine manifolds $U,V$, the Thom-Sebastiani sum $f\oplus g:U\times V\rightarrow \CC$ is also M-tame (\resp.\ cohomologically tame).
\end{prop}

\begin{rem}
Let $U=\CC_{}^{n+1}$ or $U=(\CC^*)_{}^{n+1}$. A (Laurent) polynomial $f:U\rightarrow \CC$ is {\em tame} (\cite{Broughton88}) if for some compact $K\subset U$ there exists $\varepsilon>0$ such that $\norme{\partial f(x)}\geq \varepsilon$ for $x\not\in K$. Such functions are cohomologically tame. Tameness is easily seen to be stable by Thom-Sebastiani sums, so the proposition is clear for tame (Laurent) polynomials. 
\end{rem}

\begin{proof}[Proof for M-tame functions]
Let $U$ be embedded in $\CC^M$ as a closed submanifold. Recall that $f:U\rightarrow \CC$ is not M-tame if there exists $x^o\in \PP^M-U$ and a real analytic path $x:[0,\varepsilon[\rightarrow X$ with $x(0)=x^o$, $x\lefpar]0,\varepsilon[\rigpar\subset U$, $\lim_{s\rightarrow 0}^{}f(x(s))$ exists in $\CC$ and such that for any $s\in{}]0,\varepsilon[$ one has $\ker df(x(s))\subset\ker d'\delta(x(s))$ if $d'$ denotes the holomorphic part of the differential.

Assume that $f,g$ are M-tame but not $h=f+ g$. Choose a path $z(s)=(x(s),y(s))$ as above. Denote $\delta$ the square of the distance function on $\CC^M$, $\CC^N$ or $\CC^{M+N}$. Then $\delta(x(s))$ does not remain bounded: otherwise $f(x(s))$ would remain so, as well as $g(y(s))$ because $\lim f(x(s))+g(y(s))$ exists in $\CC$; but $y(s)$ does not remain bounded, and along $y(s)$ one has $\ker dg(x(s))\subset\ker d'\delta(y(s))$; this contradicts the M-tameness of $g$. Analogously $\delta(y(s))$ does not remain bounded.

Put $x(s)=x_0s^\alpha+{}$ higher order terms and $y(s)=y_0s^\beta+{}$ higher order terms, with $x_0,y_0\neq0$. Then it follows that $\alpha,\beta<0$. Now $\dpl\frac{dh(z(s))}{ds}=dh(z(s))\cdot z'(s)$ and the inclusion $\ker dh(z(s))\subset\ker d'\delta(z(s))$ shows that there exists $\lambda(s)$ with $dh(z(s))\cdot\xi=\lambda(s)\sum_i\ov z_i(s)\xi_i$ for all $\xi$ tangent to $U\times V$. Choose $\xi=z'(s)$ and deduce that $\dpl\frac{dh(z(s))}{ds}=\lambda(s)(\alpha\norme{x_0}^2s_{}^{2\alpha-1}+\beta\norme{y_0}^2s_{}^{2\beta-1})+{}$ higher order terms.

By assumption $h(x(s))$ has a finite limit. The same is then clearly true for $\dpl\frac{dh(z(s))}{ds}$. Assume for instance that $\alpha\leq \beta$. The order of $\dpl\frac{dh(z(s))}{ds}$ is therefore equal to ${\rm ord}(\lambda)+2\alpha-1$, \ie.\ to the order of $\dpl\frac{df(x(s))}{ds}$, and this implies that $\dpl\frac{df(x(s))}{ds}$ has a finite limit, as well as $f(x(s))$. This contradicts the M-tameness of $f$.
\end{proof}

\begin{proof}[Proof for cohomologically tame functions]
Let $F:\ov X\rightarrow \PP^1$ and $G:\ov Y\rightarrow \PP^1$ be compactifications of $f,g$, such that $f$ (\resp.\ $g$) is cohomologically tame with respect to $X\defin F_{}^{-1}(\CC)$ (\resp.\ $Y\defin G_{}^{-1}(\CC)$), and let $Z\subset \ov X\times \ov Y\times \CC$ be the closure of $\{(x,y,w)\in U\times V\times \CC \mid f(x)+g(y)=w\}$. Then the map $H:Z\rightarrow \CC$ induced by the third projection is a quasi-projective compactification of $f\oplus g$. Let us show that $f\oplus g$ is cohomologically tame with respect to $Z$.

\medskip
Consider first the restriction of $H$ to the inverse image of $\CC\times \CC$ by the projection $Z\rightarrow \ov X\times \ov Y\rightarrow \PP^1\times \PP^1$, namely the map $F\oplus G:X\times Y\rightarrow \CC$. Embed $F:X\rightarrow \CC$ in $F':X'\rightarrow \CC$ with $X'$ smooth quasi-projective and $X$ closed in $X'$, and the same for $G$. Let $\ccF,\ccG$ be the direct image sheaves $j'_!\QQ_U,j'_!\QQ_V$ where now $j':U\hookrightarrow X',V\hookrightarrow Y'$ denotes the locally closed immersion.

Cohomological tameness of $f$ also means that (see \cite[\T8]{Bibi96b}), if we consider the map $T^*F':{F'}^*T^*\CC\rightarrow T^*X'$, the set $T^*{F'}_{}^{-1}(\car\ccF)\subset {F'}^*T^*\CC$ is contained in the zero section when restricted to the complement of the critical set $\crit f$ of $f$ on $U$, where $\car\ccF$ denotes the characteristic variety (or microsupport, see \eg.\ \cite{K-S90}) of $\ccF$.

From the inclusion $\car j'_!\QQ_{U\times V}^{}\subset\car\ccF\times \car\ccG$ (\cite[prop\ptbl 5.4.1]{K-S90}) we deduce that the space $T^*(F'\oplus G')_{}^{-1}(\car j'_!\QQ_{U\times V}^{})$ is contained in the zero section $T^*_\CC\CC$ when restricted to $U\times V-\crit(f\oplus g)=U\times V-\crit f\times \crit g$.

\medskip
Remark now that the projection of $Z$ in $\PP^1\times \PP^1$ does not cut $\{\infty\}\times \CC\cup\CC\times \{\infty\}$, so we may restrict to small discs $D,\Delta$ near $\infty$ in $\PP^1$ and consider $H:Z_{|D\times \Delta}^{}\rightarrow \CC$. Put $F=1/\varphi$, $G=1/\psi$ in a local coordinate on $D,\Delta$, so that $Z_{|D\times \Delta}^{}$ is defined by $\varphi(x)+\psi(y)=w\varphi(x)\psi(y)$. Let us fix $w^o\in\CC$. We will choose $D$ small enough so that $1-w\varphi(x)$ is invertible for $w$ in a neighbourhood of $w^o$ and $\varphi(x)\in D$.

Let $\wt X\rightarrow \ov X$ be a resolution of singularities such that $\wt X-U$ is a divisor with normal crossings and let $\wt Z$ be defined as $Z$ in $\wt X\times \ov Y\times \CC$. The natural map $\wt Z\rightarrow Z$ is then proper. Let $\wt\jmath$ denote the inclusion $U\times V\hookrightarrow \wt Z$. It is then enough to show that $\wt\jmath_!\QQ_{U\times V}^{}$ has no vanishing cycles with respect to the function $\wt H:\wt Z_{|D\times \Delta}^{}\rightarrow \CC$ at any point of $\wt Z$ over $(\infty,\infty,w^o)\in D\times \Delta\times \CC$: indeed, this will imply that $j_!\QQ_{U\times V}^{}$ has no vanishing cycles with respect to $H$ at any point of $Z$ over $(\infty,\infty,w^o)$, as the vanishing cycle functor commutes with the direct image functor by a proper morphism.

If $(x^o,y^o,w^o)$ is such a point, take coordinates $(x,x',x'')$ on $\wt X$ near $x^o$ such that $\wt \varphi(x,x',x'')=x^m$, where $m$ is a multi-index, and $U$ is defined by the non vanishing of the two sets of variables $x,x'$. Then $\wt Z_{|D\times \Delta}^{}$ is locally defined by $\psi(y)=x^m/(1-wx^m)$, so the map $\wt H:\wt Z_{|D\times \Delta}^{}\rightarrow \CC$ (projection on the $w$ variable) is locally analytically trivial by an isomorphism compatible with the decomposition induced by the natural stratification of $\wt X$ and the decomposition $V,\ov Y-V$ of $\ov Y$ (take the new coordinate $\wt x_1$ such that $\wt x_{1}^{m_1}=x_{1}^{m_1}/(1-wx^m)$). Consequently, $\wt\jmath_!\QQ_{U\times V}^{}$ has no vanishing cycles with respect to $\wt H$ at such a point.
\end{proof}

Let us now consider the spectrum. Recall that if $\nu',\nu'':\QQ\rightarrow \NN$ are two integral valued functions, we put $(\nu'\star\nu'')_\beta=\sum_{\beta'+\beta''=\beta}^{}\nu'_{\beta'}\nu''_{\beta''}$.

\begin{prop}
Assume that $f,g,f\oplus g$ are weakly tame. Then the spectrum at infinity of $f\oplus g$ is given by $$\nu(f\oplus g)=\nu(f)\star \nu(g).$$
\end{prop}

\begin{proof}
According to \cite[prop\ptbl 3.7]{Bibi96b} and the symmetry of the functions $\nu(f)$ and $\nu(g)$, the relation $\nu(f\oplus g)=\nu(f)\star\nu(g)$ is a consequence of the lemma below.
\end{proof}

We use notation of \loccit.:

\begin{lem}
Assume that $f,g,f\oplus g$ are weakly tame. Then the Brieskorn lattices satisfy $$G_0(f\oplus g)\simeq G_0(f)\otimes_{\CC[\theta]}^{} G_0(g).$$
\end{lem}

\begin{proof}
Recall that the complex $\lefpar \Omega_{}^{\bbullet+\dim U}(U)[\theta,\theta_{}^{-1}]e_{}^{-f/\theta},d\rigpar$ has nonzero cohomology in degree $0$ at most, as $f$ is weakly tame, and the Gauss-Manin system $G(f)$ of $f$ is by definition this cohomology space. We clearly have $G(f\oplus g)=G(f)\otimes_{\CC[\theta,\theta_{}^{-1}]}^{}G(g)$, as this is true at the level of complexes. The isomorphism of $\CC[\theta]$-modules
\begin{eqnarray*}
\Omega_{}^{\dim U}(U)[\theta]e_{}^{-f/\theta}\ootimes_{\CC[\theta]}^{} \Omega_{}^{\dim V}(V)[\theta]e_{}^{-g/\theta}&\isom&\Omega_{}^{\dim U\times V}(U\times V)[\theta]e_{}^{-(f\oplus g)/\theta}
\end{eqnarray*}
induces a surjective morphism of $\CC[\theta]$-modules $G_0(f)\otimes_{\CC[\theta]}^{} G_0(g)\rightarrow G_0(f\oplus g)$ and, as both terms are $\CC[\theta]$-free of rank $\mu(f)\mu(g)$ (weak tameness of $f,g,f\oplus g$), this morphism is an isomorphism (recall that $G_0(f)$ denotes the image of $\Omega_{}^{\dim U}(U)[\theta]e_{}^{-f/\theta}$ in $G(f)$ and that it is a free $\CC[\theta]$-module of rank $\mu(f)$, see \cite{Bibi96b}).
\end{proof}

Let $k\geq 2$ and $f:U\rightarrow \CC$ be weakly tame. Let $\sigma:\CC\rightarrow \CC$ be the multiplication by $\exp(2i\pi/k)$. Then $\sigma$ induces an automorphism $\sigma^*:G_0(f\oplus z^k)\rightarrow G_0(f\oplus z^k)$ of finite order $k$ and we put $G_0(f\oplus z^k)_{}^{(\ell)}=\ker (\sigma^*-e_{}^{2i\ell\pi/k}\id)$. For $\beta\in\QQ$ we define $\nu_\beta(G_0(f\oplus z^k)_{}^{(\ell)})$ in an evident way, as well as $\Sigma_\gamma(G_0(f\oplus z^k)_{}^{(\ell)})$ for $\gamma\in\QQ$. From the lemma we get:

\begin{cor}[\cite{Varchenko83b}]
We have $G_0(f\oplus z^k)=\oplus_{\ell=0}^{k-1}G_0(f\oplus z^k)_{}^{(\ell)}$. For any $\beta\in\QQ$ we have $$\nu_\beta(G_0(f\oplus z^k)_{}^{(\ell)})=\nu_{\beta-\ell/k}^{}(f)$$ and for any $\gamma\in\QQ$ we have $$\Sigma_\gamma(G_0(f\oplus z^k)_{}^{(\ell)})=\Sigma_{\gamma-\ell/k}^{}(f).$$ \qed
\end{cor}

\section{Positivity}\label{secpos}
Let $f:U\rightarrow \CC$ be a weakly tame function. We will consider the sheaves on $\CC$ with fiber at $c\in\CC$ equal to $H^i(U,f_{}^{-1}(c);\QQ)$. More precisely, let $k:W_f\hookrightarrow U\times \CC$ be the open embedding complementary to $\graph (f)\subset U\times \CC$ and let $p:U\times \CC\rightarrow \CC$ be the projection. Let $\ccF=\bR p_* k_!\pQQ_{W_f}$ where $\pQQ_{W_f}=\QQ_{W_f}[\dim {W_f}]$.

\begin{lem}
For $c\in\CC$ we have $\cH^i\ccF_c=H^{i+\dim U+1}(U,f_{}^{-1}(c);\QQ)$ and this space is equal to $0$ if $i\neq -1$, hence $\ccF$ has cohomology in degree $-1$ only. Moreover $\ccF$ has perverse cohomology in degree $0$ only.
\end{lem}

\begin{proof}
Consider the triangle
$$
\bR f_*\pQQ_U\MRE{}\bR p_* k_!\pQQ_{W_f}\MRE{} \bR p_*\pQQ_{U\times \CC}^{}\MRE{+1}.
$$
Use that $\bR f_*$ commutes with the base change by $\{c\}\hookrightarrow \CC$ (see \cite[cor\ptbl 8.4]{Bibi96b}) and that $\bR p_*$ does so for the third term and conclude that $\bR p_*$ does so for the second term to get the first statement. The fibre at $c$ of $\cH^i\ccF$ is thus $H^{i+\dim U+1}(U,f_{}^{-1}(c);\QQ)$. Arguments analogous to the ones of \cite[lemma 8.5]{Bibi96b} give that $\cH^i\ccF=0$ for $i\geq 0$ and are constant sheaves for $i<-1$. Using Leray spectral sequence and the fact that $H^*(U\times \CC,\graph (f);\QQ)=0$ one gets $\cH^i\ccF=0$ for $i<-1$. In order to get the perversity statement, it is enough to show that for each $c\in\CC$ the complexes $\psip_{t-c}^{}\ccF$ and $\phip_{t-c}^{}\ccF$ have cohomology in degree $0$ only. As $\ccF$ is, generically on $\CC$, a perverse sheaf (\ie.\ a locally constant sheaf up to a shift by one), the result is true for $\psip$. On the other hand we have $\phip_{t-c}^{}\ccF=\phip_{t-c}^{}\bR f_*\pQQ_U$ and weak tameness implies that this complex has cohomology in degree $0$ only (see \cite[cor\ptbl 8.4]{Bibi96b}).
\end{proof}

Consider now a family $f:\cU\rightarrow \CC$, with $\pi:\cU\rightarrow S$ a smooth affine morphism over a germ of curve $S$ with coordinate $s$, and denote $k:{\cW_f}\hookrightarrow \cU\times \CC$ the complementary inclusion of the graph of $f$. We still denote $\ccF=\bR (p,\pi)_*k_!\pQQ_{\cW_f}$ which is now a complex on $\CC\times S$.

\begin{lem}
Assume that $f_0:\cU_0\rightarrow \CC$ is weakly tame. Then for $m<0$ the variation mapping $\phip_s\per\cH^m\ccF\rightarrow \psip_s\per\cH^m\ccF$ is an isomorphism and for $m=0$ we have an exact sequence
$$
0\longrightarrow \phip_s\per\cH^0\ccF\longrightarrow \psip_s\per\cH^0\ccF\longrightarrow \per\cH^0\per i_0^!\ccF\longrightarrow 0
$$
where $i_0$ denotes the inclusion $\{0\}\hookrightarrow S$.
\end{lem}

\begin{proof}
Recall that we denote here, for a complex $\ccG$,  $\per i_{}^{-1}\ccG=i_{}^{-1}\ccG[-1]$ and $\per i^!\ccG=i^!\ccG[1]$, where $i_{}^{-1}$ and $i^!$ are the standard functors of sheaf theory. As $\pi:\cU\rightarrow S$ is smooth, we have $\per i_{\pi_{}^{-1}(0)}^{!}\pQQ_\cU=\per i_{\pi_{}^{-1}(0)}^{-1}\pQQ_\cU=\pQQ_{\cU_0}$ and $\per i_{\pi_{}^{-1}(0)}^{!}\pQQ_{\cU\times \CC}=\per i_{\pi_{}^{-1}(0)}^{-1}\pQQ_{\cU\times \CC}=\pQQ_{\cU_0\times \CC}$, so
$$
\per i_{\pi_{}^{-1}(0)}^{!}k_!\pQQ_{\cW_f}=\per i_{\pi_{}^{-1}(0)}^{-1}k_!\pQQ_{\cW_f}= k_!\per i_{\pi_{}^{-1}(0)}^{-1}\pQQ_{\cW_f}=k_!\pQQ_{W_{f_0}}^{}
$$
and $\per i_0^!\bR(p,\pi)_*k_!\pQQ_{\cW_f}=\bR p_*k_!\pQQ_{W_{f_0}}^{}$ is perverse, according to the previous lemma. The lemma follows now from the long exact sequence in perverse cohomology associated with the triangle
$$
\phip_s\ccF\MRE{}\psip_s\ccF\MRE{}\per i_0^!\ccF\MRE{+1}
$$
the fact that $\per\cH^m\ccF=0$ for $m>0$ and the fact that $\phip_s$ and $\psip_s$ commute with taking perverse cohomology (see \cite{Brylinski86}). 
\end{proof}

Keep notation as above and for $c\in\CC$ let $\ccF_c$ be the inverse image of $\ccF$ by the inclusion $\{c\}\times S\hookrightarrow \CC\times S$ shifted by $-1$. Then for $c$ general enough this inclusion is noncharacteristic for $\ccF$ and $\per\cH^m(\ccF_c)=(\per\cH^m\ccF)_c$. For such a $c$ we have $\ccF_c=\bR \pi_*k_!\pQQ_{\cU-\{f=c\}}^{}$; the cohomology of the fibre at $s\neq0$ of this complex is the relative cohomology $H^\bbullet(\cU_s,f_{s}^{-1}(c);\QQ)$. If $f_0$ is weakly tame, the results of the previous lemma hold for $\ccF_c$ so that $\per i_0^!\ccF_c$ has cohomology in degree $0$ only, this cohomology being the relative cohomology $H_{}^{\dim \cU_0}(\cU_0,f_{0}^{-1}(c);\QQ)$.

\begin{poslem}
Assume that for any $s\in S$ the function $f_s:\cU_s\rightarrow \CC$ is weakly tame. Then for $c\in\CC$ general enough, $\ccF_c$ is a perverse sheaf on $S$ and we have an exact sequence of vector spaces
$$
0\longrightarrow \phip_s\ccF_c\longrightarrow \psip_s\ccF_c\longrightarrow \per i_0^!\ccF_c\longrightarrow 0.
$$
\end{poslem}

\begin{proof}
By noncharacteristic restriction, we may replace $\ccF$ with $\ccF_c$ in the previous lemma. As $f_s$ is weakly tame for any $s\in S_\{0\}$, the perverse sheaf $\per\cH^m\ccF_c$ is supported at $s=0$ for $m<0$. Hence $\psip_s\per\cH^m\ccF_c=0$, and by the previous lemma we also have $\phip_s\per\cH^m\ccF_c=0$, so $\per\cH^m\ccF_c=0$ for $m<0$ and $\ccF_c=\per\cH^0\ccF_c$. The previous lemma then shows that $\per i_0^!\ccF_c$ is also perverse and gives the desired exact sequence.
\end{proof}

\section{Proof of the theorem}\label{secproofth}
Let $V$ be any smooth quasi-projective manifold of pure dimension $\dim V$, let $Z$ be a closed subvariety and let $F_{D}^{\bbullet}H^k(V,Z;\CC)$ be the (decreasing) Hodge-Deligne filtration on the cohomology spaces of the pair $(V,Z)$. For $q\in\NN$, let
\begin{eqnarray*}
\chi_{\rm Del}^{}(V,Z;q)&\defin&\sum_i(-1)^i\dim\gr_{F_D}^{q}H^{\dim V-i}(V,Z;\CC)\\
\zeta_{\rm Del}^{}(V,Z;S)&\defin&\prod_q(S+\dim V-q)_{}^{\chi_{\rm Del}^{}(V,Z;q)}.
\end{eqnarray*}
If $Z$ is a hypersurface of $V$, we have $\chi_{\rm Del}^{}(V,Z;q)=\chi_{\rm Del}^{}(V;q)+\chi_{\rm Del}^{}(Z;q)$.

If $(V,Z)$ comes equipped with an automorphism $\sigma$ of finite order $k$, the Hodge filtration on the cohomology splits with respects to the eigenvalues $e_{}^{2i\pi\ell/k}$ of $\sigma^*$ and we can define $\chi_{\rm Del}^{(\ell)}(V,Z;q)$ and $\zeta_{\rm Del}^{(\ell)}(V,Z;S)$ for each $\ell$ between $0$ and $k-1$.

For any $R\in\CC(S)$ with zeros and poles in $\QQ$ (or in $\RR$) and for $\gamma\in\QQ$, put

$\bullet$ $\deg_{\gamma}^{\psi}R(S)$: sum of the degrees of the terms $(S+\beta)$ in $R$ with $\beta\in[\gamma,\gamma+1[$,

$\bullet$ $\deg_{\gamma}^{\phi}R(S)$: sum of the degrees of the terms $(S+\beta)$ in $R$ with $\beta\in{}]\gamma,\gamma+1]$.

\medskip
If $R$ has only integral poles or zeros, $\deg_{\gamma}^{\psi}R(S)$ (\resp.\ $\deg_{\gamma}^{\phi}R(S)$) is the degree of $S+\gamma$ (\resp.\ $S+\gamma+1$) in $R$, also denoted $\deg_{S+\gamma}^{}R$ (\resp.\ $\deg_{S+\gamma+1}^{}R$).

\medskip
In the following we use notation of \cite[\T5]{Bibi96a}.
\begin{prop}
Let $f:\cU\rightarrow \CC$ be a one parameter family of regular functions on a smooth affine pure dimensional family $\pi:\cU\rightarrow S$, such that each $f_s:\cU_s\rightarrow \CC$ is weakly tame. Then we have for any $\gamma\in\ZZ$, $k\geq 2$, $0\leq \ell\leq k-1$ and $c\in\CC$ general enough, putting $q=\dim\cU_s-1-\gamma$,
\begin{eqnarray*}
\Sigma_{\gamma-\ell/k}^{}(f_s)-\Sigma_{\gamma-\ell/k}^{}(f_0)&=&\chi_{\rm Del}^{(\ell)}(\cU_s\times \CC,(f_{s}\oplus z^k)^{-1}(c);q)- \chi_{\rm Del}^{(\ell)}(\cU_0\times \CC,(f_{0}\oplus z^k)^{-1}(c);q).
\end{eqnarray*}
\end{prop}

\begin{proof}
Remark first that the LHS can be written $\dpl\deg_{\gamma}^{\phi}\lefpar \frac{\SP_\psi(G_{0}(f_{s}\oplus z^k)^{(\ell)};S)}{\SP_\psi(G_{0}(f_{0}\oplus z^k)^{(\ell)};S)}\rigpar$, as follows from the last corollary in \T\ref{secTS}. 

Let us show first that, for $q=\dim \cU_0-1-\gamma$, we have
\begin{equation}\tag{$*$}
\deg_{\gamma}^{\phi}\lefpar \frac{\SP_\psi(G_{0}(f_{s})^{};S)}{\SP_\psi(G_{0}(f_{0})^{};S)}\rigpar=\chi_{\rm Del}^{}(\cU_s,f_{s}^{-1}(c);q)- \chi_{\rm Del}^{}(\cU_0,f_{0}^{-1}(c);q).
\end{equation}
If $f:U\rightarrow \CC$ is a weakly tame regular function on a smooth affine manifold, we have (\cite[cor\ptbl 13.6]{Bibi96b})
\begin{equation}\tag{$**$}
\SP_\psi(G,G_0;S)=\SP_\phi(\psip_{1/t}\bR f_*\pCC_U;S)\cdot \zeta_{\rm Del}^{}(U;S)
\end{equation}
and taking degrees we get, for $\gamma\in\QQ$
\begin{eqnarray*}
\Sigma_\gamma(f)&=&\deg_{\gamma}^{\phi}\SP_\psi(G,G_0;S)\\
&=&\deg_{\gamma}^{\phi}\SP_\phi(\psip_{1/t}\bR f_*\pCC_U;S)+ \deg_{\gamma}^{\phi} \zeta_{\rm Del}^{}(U;S)\\
&=&\deg_{\gamma}^{\psi}\SP_\psi(\psip_{1/t}\bR f_*\pCC_U;S)+ \deg_{\gamma}^{\phi} \zeta_{\rm Del}^{}(U;S).
\end{eqnarray*}

Now, from the theory of mixed Hodge modules \cite{MSaito87} (or the property of the limit mixed Hodge structure at infinity, as constructed by Steenbrink and Zucker \cite{S-Z85}) we know that for $\gamma\in\ZZ$ we have
\begin{eqnarray*}
\deg_{\gamma}^{\psi}\SP_\psi(\psip_{1/t}\bR f_*\pCC_U;S)&=&\deg_{S+\gamma}\zeta_{\rm Del}^{}(f_{}^{-1}(c);S)
\end{eqnarray*}
for $c$ general enough. Equality $(*)$ follows easily.

\medskip
We may now argue as in \cite{Varchenko83}: we may apply the previous  argument to $f\oplus z^k$ which is weakly tame and we need to show that all identifications made to get $(**)$ remain true on the $e_{}^{2i\pi\ell/k}$-eigenspaces. We will indicate the main steps, following the proofs in \cite{Bibi96a,Bibi96b}.

Put $g=f\oplus z^k:V=U\times \CC\rightarrow \CC$ and choose a projective compactification $\cY$ of $V$ on which $g$ extends as $G:\cY\rightarrow \PP^1$ and $\sigma$ extends as an automorphism $\sigma$ which commutes with $G$. So $\sigma$ induces an automorphism $\sigma$ of $Y=G_{}^{-1}(\CC)$ wich commutes with $G$ and with $j:V\hookrightarrow Y$. Using an equivariant version of Hironaka's resolution theorem of singularities, we may also assume that $\cY-V$ is a divisor with normal crossings.

We have a natural isomorphism $\QQ_V\simeq \sigma_*\QQ_V$. We deduce from it a lifting $\lambda^\sigma:\bR j_*\pQQ_V\isom\sigma_*\bR j_*\pQQ_V$ in the category $\Perv (Y)$ of $\QQ$-perverse sheaves on $Y$, which is filtered with respect to the weight filtration $W_\bbullet\bR j_*\pQQ_V$.

Let $\sigma$ still denote the automorphism $\sigma\times \id$ of $\cY\times \Afuc$, where $\Afuc$ denotes the affine line whith coordinate $\tau$ (the Fourier plane). We introduced in \cite[\T1.5]{Bibi96a} functors $\cFG,\FGir:\Perv (Y)\rightarrow \Perv (\cY\times \Afuc)$. One verifies that both functors commute with $\sigma_*$.

In an analogous way one verifies that $\sigma_*$ commutes with the nearby and vanishing cycles functors $\psip_{1/g}^{}$, $\psip_{\tau,1}^{}$, $\phip_\tau$ and the corresponding monodromies.

The isomorphism $\phip_\tau\circ\cFG\simeq \psip_{1/t}$ of \cite[cor\ptbl 1.13]{Bibi96a} being functorial, we get a commutative diagram
$$
\begin{array}{ccc}
\phip_\tau\cFG\bR j_*\pQQ_V&\simeq&\psip_{1/G}\bR j_*\pQQ_V\\
\MDL{\phip_\tau\cFG(\lambda^\sigma)}&&\MDR{\psip_{1/G}(\lambda^\sigma)}\\
\phip_\tau\cFG\sigma_*\bR j_*\pQQ_V&\simeq&\psip_{1/G}\sigma_*\bR j_*\pQQ_V
\end{array}
$$
and the proof of \cite[cor\ptbl 1.13]{Bibi96a} shows that the isomorphism $\phip_\tau\circ\cFG\simeq \psip_{1/t}$ is compatible with the one expressing the commutation with $\sigma_*$, because both act on different sets of variables. We finally get a commutative diagram of liftings
$$
\begin{array}{ccc}
\phip_\tau\cFG\bR j_*\pQQ_V&\simeq&\psip_{1/G}\bR j_*\pQQ_V\\
\MDL{\lambda^\sigma_\tau}&&\MDR{\lambda^\sigma_{1/G}}\\
\sigma_*\phip_\tau\cFG\bR j_*\pQQ_V&\simeq&\sigma_*\psip_{1/G}\bR j_*\pQQ_V
\end{array}
$$
and the fact that these liftings are morphisms of mixed Hodge modules can be obtained by a local computation on $\cY$. As indicated above, this is sufficient to give the proposition.
\end{proof}

\begin{proof}[End of the proof of the theorem]
Apply the positivity lemma of \T\ref{secpos} to the family $f\oplus z^k$ (recall from \T\ref{secTS} that $f_s\oplus z^k$ is weakly tame for any $s\in S$). The surjective morphism given by the lemma splits with respect to the eigenspaces of $\lambda^\sigma$ and is strict with respect to the Hodge filtrations as it is induced by a morphism of mixed Hodge modules (see \cite{MSaito87}). Using once more that the Hodge numbers of $\psip_s$ are equal to the ones of the general fibre at $s\neq0$, we get that for any $q\in\ZZ$ we have
$$
\chi_{\rm Del}^{(\ell)}(\cU_s\times \CC,(f_{s}\oplus z^k)^{-1}(c);q)\geq \chi_{\rm Del}^{(\ell)}(\cU_0\times \CC,(f_{0}\oplus z^k)^{-1}(c);q)
$$
so the RHS in the proposition is $\geq 0$.
\end{proof}

\nocite{Mebkhout87}

\providecommand{\bysame}{\leavevmode\hbox to3em{\hrulefill}\thinspace}

\end{document}